\documentclass[12pt]{article}
\usepackage[T2A]{fontenc}
\usepackage[cp1251]{inputenc}
\usepackage[english]{babel}
\usepackage{amssymb,amsmath,amsthm,amsfonts}
\usepackage{color}
\usepackage{cite}
\usepackage[a4paper,left=20mm,right=20mm,bottom=20mm,top=20mm]{geometry}
\usepackage{colordvi}

\theoremstyle{plain}

\theoremstyle{definition}

\begin{document}

\bigskip


\title{A note on the supersolvability of a finite group with prime index of some subgroups}
\author{V.\,S.~Monakhov, A.\,A.~Trofimuk}

\date{January 16, 2019}

\maketitle

\begin{abstract}
In this paper, we proved that a group $G$ is supersoluble if and only if for any prime $p\in \pi (G)$ there exists a supersoluble subgroup of
index $p$.
\end{abstract}

All groups considered in this paper will be finite.  Denote by $\pi(G)$~the set of all prime divisors of order of~$G$. The notation $H\le G$ means that~$H$~is a subgroup of~$G$.

By Huppert \cite[Theorem~6]{hup54}, a group $G$ is supersoluble if and only if every maximal subgroup of $G$ has prime index. The following observation is generated by Huppert's result.

{\bf Theorem~1.} {\sl
Let  $G$ be a group. Then $G$ is supersoluble if and only if for any prime $p\in \pi(G)$ there exists a supersoluble subgroup of index~$p$. }

{\bf Proof.}
If $G$ is supersoluble, then for every  $p\in \pi(G)$ maximal subgroup of $G$ that contains a Hall $p^{\prime}$-subgroup $G_{p^{\prime}}$ is supersoluble and has index~$p$.

Conversely, let $\pi(G)=\{p_1, p_2, \ldots, p_n\}$, $p_1>p_2>\ldots >p_n$ and $H_i$ be a supersoluble subgroup such that $|G:H_i|=p_i$, where $i=1,2,\ldots, n$.
It is clear that $H_n$ is normal in $G$.  Let $P$ be a Sylow  $p_1$-subgroup of~$H_n$. Then $P$ is normal in $G$. We use induction on the order of $G$. Suppose that $\Phi(P)\neq 1$. Then
$$
\Phi(P)\leq \Phi(G)\leq \bigcap_{\substack{i=1,2,\ldots, n}} H_i
$$
and $G/\Phi(G)$ is supersoluble by induction. So $G$ is supersoluble.
Therefore, we consider that $\Phi(P)= 1$ and $P$ is elementary abelian.
By Mashke's theorem \cite[Theorem~I.17.6]{Hup},  $P=N_1\times N_2\times \ldots\times N_m $, where $N_j$ is a minimal normal subgroup of~$G$, where $i=1,2,\ldots, m$.
Since $|G:H_1|=p_1$, it follows that $P$  is not contained in $H_1$ and there is a subgroup $N_j$ for some $j$ such that it is not contained in~$H_1$. Hence $G=N_jH_1$,
$N_j\cap H_1=1$ and $|N_j|=p_1$. Now $G$ is supersoluble by \cite[Lemma~4.46]{Mon}. The theorem is proved.

{\bf Corollary~1.} (\cite[Theorem~3.1]{As}) {\sl
Let~$H$ and~$K$~be supersoluble subgroups of~$G$ and $G=HK$.
If each subgroup of~$H$ permutes with every subgroup of~$K$, then~$G$ is
supersoluble.}

{\bf Proof.} We use induction on the order of $G$.
We can assume that $H_1K\ne G\ne HK_1$ for all proper subgroups~$H_1$ of~$H$
and~$K_1$ of~$K$. Let $p\in \pi (H)$, $A\le H$ and $|H:A|=p$.
By induction,~$AK$ is supersoluble and~$|G:AK|=p$.
Similarly, if $q\in \pi (K)$, $B\le K$ and $|K:B|=q$,
then~$HB$ is supersoluble and~$|G:HB|=q$.
Since $\pi (H)\cup \pi (K)=\pi (G)$, it follows that~$G$ is supersoluble
by Theorem~1.

\medskip

Asaad and Shaalan's theorem \cite[Theorem~3.1]{As} was developed by Guo~W., Shum ~K.\,P., Skiba~A.\,N., see~\cite[Theorem~A]{Skiba1}. This result can also be obtained from Theorem~1. We introduce the following

\medskip
{\bf Definition.}
The subgroups $H$ and $K$ are said to be  {\sl $tcc$-permutable} if for any $X\le H$ and $Y\le K$ there exists an element $u\in \langle X,Y\rangle $ such that~$XY^u=Y^uX$.
Note that the equality~$ XY^u=Y^uX $ is equivalent to~$XY^u$~is a subgroup.

{\bf Lemma~1.} {\sl
Let~$H$ and~$K$~ be subgroups of~$G$.
If $H$ and $K$ are  $tcc$-permutable, then the subgroups $H$ and~$K^h$ are
$tcc$-permutable for any $h\in H$.
}

{\bf Proof.}
Let $X\le H$ and $Y\le K^h$. Then $X^{h^{-1}}\le H$, $Y^{h^{-1}}\le K$
and there exists an element $u\in \langle X^{h^{-1}},Y^{h^{-1}}\rangle $ such that
$X^{h^{-1}}(Y^{h^{-1}})^u$~is a subgroup.
Since $\langle X^{h^{-1}},Y^{h^{-1}}\rangle =\langle X,Y\rangle ^{h^{-1}}$,
it follows that there exists an element  $v\in \langle X,Y\rangle $ such that $u=v^{h^{-1}}$.
Because $h^{-1}u=vh^{-1}$, we have
$$
X^{h^{-1}}(Y^{h^{-1}})^u=X^{h^{-1}}Y^{vh^{-1}}
=(XY^v)^{h^{-1}}.
$$
Hence $XY^v$ is a subgroup. Consequently, $H$ and $K^h$ are $tcc$-permutable.

\medskip

{\bf Corollary~2.} (\cite[Theorem~A]{Skiba1}) {\sl
Let~$H$ and~$K$~be supersoluble subgroups of~$G$ and~$G=HK$. If for any~$H_1\leq H$ and $K_1\leq K$ there exists an element $x\in \langle H_1,K_1\rangle $ such that $H_1K_1^x=K_1^xH_1$,
then~$G$ is supersoluble.}

{\bf Proof.}
It is obvious that~$H$ and~$K$ are $tcc$-permutable.
We use induction on the order of $G$. We can assume that $H_1K\neq G\neq HK_1$ for any proper subgroups~$H_1$ of~$H$ and~$K_1$ of~$K$.

Let $p\in \pi(H)$, $H_1\leq H$ and $|H:H_1|=p$. By hypothesis, there exists an element $g\in \langle H_1,K\rangle $ such that  $H_1K^g=K^gH_1$. Since $g=hk$ for some
$k\in K$ and $h\in H$, we have $H_1K^g=H_1K^h=K^hH_1$ and
$h\in \langle H_1,K\rangle $. Because $G=HK^h$, it follows that $|G:H_1K^h|=p$.
By Lemma~1, $H$ and~$K^h$ are $tcc$-permutable.
Hence the subgroups $H_1$ and $K^g$ satisfy the hypothesis of the corollary. By induction, $H_1K^g$ is supersoluble.

Similarly, if $q\in \pi(K)$, $K_1\leq K$ and $|K:K_1|=q$, then there exists an element $k\in K$ such that $H^kK_1$ is supersoluble in~$G$ and has index $q$.
Now $G$ is supersoluble by Theorem~1.

\end{document}